\documentclass{article}

\usepackage[preprint]{neurips_2020}

\usepackage[utf8]{inputenc} 
\usepackage[T1]{fontenc}    
\usepackage{hyperref}       
\usepackage{url}            
\usepackage{booktabs}       
\usepackage{amsfonts}       
\usepackage{nicefrac}       
\usepackage{microtype}      
\usepackage{graphicx,amsmath,amssymb,amsthm,listings,xcolor,caption,subcaption,bm}
\usepackage{algorithm}
\usepackage{algorithmic}
\usepackage{multirow}
\usepackage{bm}
\newcommand{\symnet}{\text{SymNet}}

\title{RODE-Net: Learning Ordinary Differential Equations with Randomness from Data}

\author{%
  Junyu Liu$^*$ \\
  Academy for Advanced Interdisciplinary Studies \\
  Peking University \\
  \texttt{liujunyu@pku.edu.cn} \\
  \And
  Zichao Long\thanks{Equal Contribution} \\
  School of Mathematical Science\\
  Peking University \\
  \texttt{zlong@pku.edu.cn} \\
  \And 
  Ranran Wang \\
  National School of Development \\
  Peking University \\
  \texttt{ranranw@bibdr.org} \\
  \And 
  Jie Sun\thanks{https://sites.google.com/view/jiesun/home} \\
  Theory Lab \\
  Huawei Hong Kong Research Center \\
  \texttt{j.sun@huawei.com} \\
  \And
  Bin Dong\thanks{http://bicmr.pku.edu.cn/~dongbin} \\
  Beijing International Center for Mathematical Research \\ Institute for Artificial Intelligence \\ Center for Data Science\\
  Peking University\\
  \texttt{dongbin@math.pku.edu.cn} \\
}

\begin{document}

\maketitle

\begin{abstract}
  Random ordinary differential equations (RODEs), i.e. ODEs with random parameters,
  are often used to model complex dynamics. Most existing methods to identify 
  unknown governing RODEs from observed data often rely on strong prior knowledge.
  Extracting the governing equations from data with less prior knowledge 
  remains a great challenge. 
  In this paper, we propose a deep neural network, called RODE-Net, to tackle 
  such challenge by fitting a symbolic expression of the differential equation 
  and the distribution of parameters simultaneously.
  To train the RODE-Net, we first estimate the parameters of the unknown RODE 
  using the symbolic networks\cite{long2019pde} by solving a set of deterministic 
  inverse problems based on the measured data, and use a generative adversarial 
  network (GAN) to estimate the true distribution of the RODE's parameters.
  Then, we use the trained GAN as a regularization to further improve the 
  estimation of the ODE's parameters. The two steps are operated alternatively. 
  Numerical results show that the proposed RODE-Net can well estimate the
  distribution of model parameters using simulated data and can make reliable 
  predictions. 
  It is worth noting that, GAN serves as a data driven regularization in RODE-Net and is more effective than the $\ell_1$ based regularization that is often used in system identifications.
\end{abstract}

\section{Introduction}
In the study of complex systems such as fluid dynamics, soft 
materials, biological evolution, etc., we often introduce empirical formulas 
and parameters into differential equations to model complex macroscopic
phenomena induced by microscopic behaviors. These include for example, the 
sub-grid stress in turbulence models, high order moment closure for moment 
equations, the drag coefficient which expressed as a function of porosity 
in particle fluid problem, and a variety of empirical formulas in biological 
and economical models. However, drawbacks of using empirical 
formulas and parameters to model such complex systems are that:
1) for systems with non-separable scales, it is difficult to directly 
deduce accurate macroscopic equations based on microscopic mechanism;
2) models with fixed parameters cannot reflect the inherent randomness 
of the system, and hence can only crudely approximate the system in 
the average sense.
Therefore, to address the above issues, we propose a machine learning
framework to learn both the expression of differential equations and 
the distribution of the associated random parameters simultaneously. 

\subsection{Related Work}
We start with a review of inverse problems for deterministic systems. 
Under the premise that the explicit form of the differential equation 
is known, the optimal parameters can be learned from the observed data 
by using classical system identification and tools from inverse problems. 
Several works \cite{lin2008learning,liu2010learning,long2018pde,
patel2018nonlinear} solved a class of inverse problems under weaker 
assumptions based on the idea of unrolling dynamics of the numerical 
integration in the time direction. In \cite{raissi2018hidden}, the authors 
proposed to use neural network as a surrogate to learn optimal parameters. 
When the explicit formula of the equation is unknown, 
\cite{brunton2016discovering} constructed a dictionary consisted of simple 
terms that are likely to appear in the equations and employ sparse 
regression methods to select candidates for the expression of the equations. 
A series of works \cite{bongard2007automated,schmidt2009distilling,xu2020dlga} used genetic 
algorithms to discover the underlying terms of the differential equations. In \cite{long2019pde}, the authors proposed
a deep symbolic neural network, called SymNet, to estimate the unknown 
expression of the differential equations, which has a relatively low memory requirement 
and computational complexity in many cases.

For problems with randomness, uncertainty 
quantification methods \cite{smith2013uncertainty} are often adopted to study forward
uncertainty propagation (e.g., \cite{liu1986probabilistic,zhang2001stochastic,
ghanem2003stochastic,xiu2010numerical}). However, it is generally much more 
difficult to study inverse problems with uncertainty than forward uncertainty 
propagation. In earlier work, \cite{kennedy2001bayesian} proposed a 4-steps 
modular Bayesian approach to calibrate parameters of models. In 
\cite{chen2019learning}, neural networks were used to estimate the modes of 
K-L expansion during the study of the forward and inverse problems of 
stochastic advection-diffusion-reaction equations. Based on PINN
\cite{raissi2018hidden}, \cite{yang2020b} estimated the posterior of model 
parameters by Hamiltonian Monte Carlo and variational inference. 
The most related work to ours is \cite{yang2018physics}, where the authors introduced a generative adversarial network (GAN) 
\cite{goodfellow2014generative,arjovsky2017wasserstein,gulrajani2017improved}
to estimate the distribution of data snapshots. These works further advanced the development of inverse uncertainty quantification. However, 
Bayesian based framework \cite{kennedy2001bayesian,chen2019learning,yang2020b} 
often suffers from the curse of dimensionality. Although the use of GAN by 
\cite{yang2018physics} can overcome the curse of dimensionality, it also requires a relatively strong prior knowledge on the differential equation to be identified. Furthermore, the GAN of \cite{yang2018physics} were not specifically designed to estimate the distribution of the parameters of the differential equations. 

In this paper, we propose a new machine learning framework for inverse uncertainty quantification. We introduce a new model, called RODE-Net, to estimate random ordinary differential equations (RODEs) from observed data by combining SymNet 
\cite{long2019pde} for system identification and GAN 
\cite{goodfellow2014generative,arjovsky2017wasserstein,gulrajani2017improved}
for parameters distribution estimation. The high expressive power of SymNet enables the proposed RODE-Net to assume only minor prior knowledge on the form of the RODE to be identified. 
A particular novelty of the propose model is that, unlike existing inverse uncertainty quantification methods, RODE-Net estimates the distribution of parameters and makes use of GAN as a data driven regularization. 

We note that, using GAN as a data driven regularization is common in image restoration
\cite{bora2017compressed,shah2018solving,ledig2017photo}.
For example, \cite{bora2017compressed,shah2018solving} put forward the idea  that we could solve image restoration problems within the range of a well trained GAN. In  \cite{ledig2017photo} the authors used discriminator as the regularization in image super resolution problems. These latest studies in computer vision inspired us to use GAN as a data driven regularization which is new to inverse uncertainty quantification.

The rest of the paper is organized as follows. In Section 2, we introduce the architecture of RODE-Net. Details on training and the loss functions are introduced in Section 3. Experiments are presented in Section 4, and we conclude the paper in Section 5.

\section{The RODE-Net}
Given a set of observed time series $\{\bm{x}(t)\} \subset \mathbb{R}^d$ 
with $d$ being the dimension of the observable quantities, we aim to 
discover the governing RODE from the set of time series. We assume that the RODE to be 
discovered takes the following form:
\begin{equation}\label{eq:RODE}
  \frac{d\bm{x}}{dt} = \bm{F}_{\bm{\eta}}(\bm{x}),
  \bm{x}=(x_1,\cdots,x_d)\in\mathbb{R}^d,t\geq0
\end{equation}
where $\bm{F}_{\bm{\eta}}(\cdot): \mathbb{R}^d\to\mathbb{R}^d$ is a $d$ dimensional 
vector function and $\bm{\eta} \sim \mathcal{P}$ is the random parameter 
of $\bm{F}_{\bm{\eta}}$. The RODE describes an infinite set of ODEs determined by the distribution of the random vector $\bm{\eta}$. For each random realization of the random variable ${\bm{\eta}}$, denoted as ${\bm{\eta}}_i$, we 
call the associated ODE as ODE-${\bm{\eta}}_i$. We note that the RODEs in this paper are different from those considered in 
\cite{han2017random,neckel2013random}. To obtain solutions with higher regularities, the latter considers system parameters as random processes (e.g. solution of SDEs) rather than random variables. 

Consider the data set $X:=\{\bm{x}_{\bm{\eta}_i}^j(t_k):1\leq i\leq M,1\leq j\leq N_i,
0\leq k\leq S\}$, where $M$ is the number of realizations of ODE-${\bm{{\bm{\eta}}}}$. In practice, the data set $X$ is measured from the real world. As a proof of concept and to validate the proposed RODE-Net, in our experiments in Section 4, we generate data $\bm{x}_{\bm{\eta}_i}^j(t_k)$ by solving ODE-${\bm{\eta}}_i$ using the fourth order Runge-Kutta method with $N_i$ different initial values $\bm{x}_{\bm{\eta}_i}^j(t_0), j=1,\cdots,N$. 

\begin{figure}[htbp!]
		\centering
		\includegraphics[width = 0.7\textwidth]{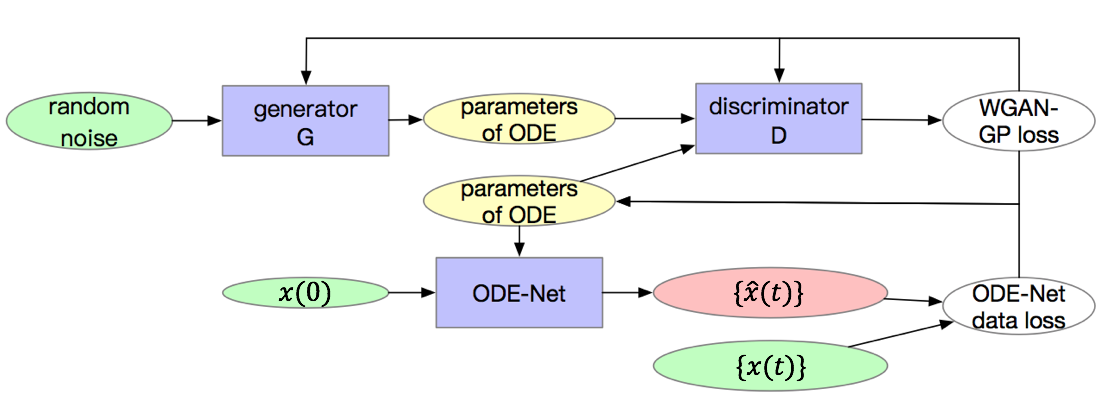}
		\caption{The pipeline of RODE-Net. The modules in green indicate input data; the modules in yellow indicate intermediate outputs; the module in red indicates the final output; the modules in blue indicate those with trainable parameters; and the modules in white indicate the loss functions.}\label{fig:RODE-net}
\end{figure}

Our proposed RODE-Net (Fig. \ref{fig:RODE-net}) is designed to identify each individual ODE of \{ODE-${\bm{\eta}}_i$ $: 1\le i\le M$\} and to estimate the distribution of ${\bm{\eta}}$ from $\{\bm{\eta}_i: 1\le i\le M\}$ at the same time. To achieve both objectives, the RODE-Net framework consists of two main components: 
\begin{enumerate}
    \item A set of SymNet-based ODE-Nets which aim to estimate $\bm{F}_{\bm{\eta}_i}$ of each individual ODE;
    \item A GAN to learn the distribution of $\{{\bm{\eta}}_i\}$.
\end{enumerate}
The trained GAN represents the distribution $\mathcal{P}$ of the RODE and is able to generate ODE instances ODE-${\bm{\eta}}$. Moreover, we note that 
the GAN also serves as a good regularization in the training of ODE-Nets.

\subsection{The ODE-Net Component}\label{subsec:ode-net}
\begin{figure}[htbp!]
		\centering
		\includegraphics[width = 0.8\textwidth]{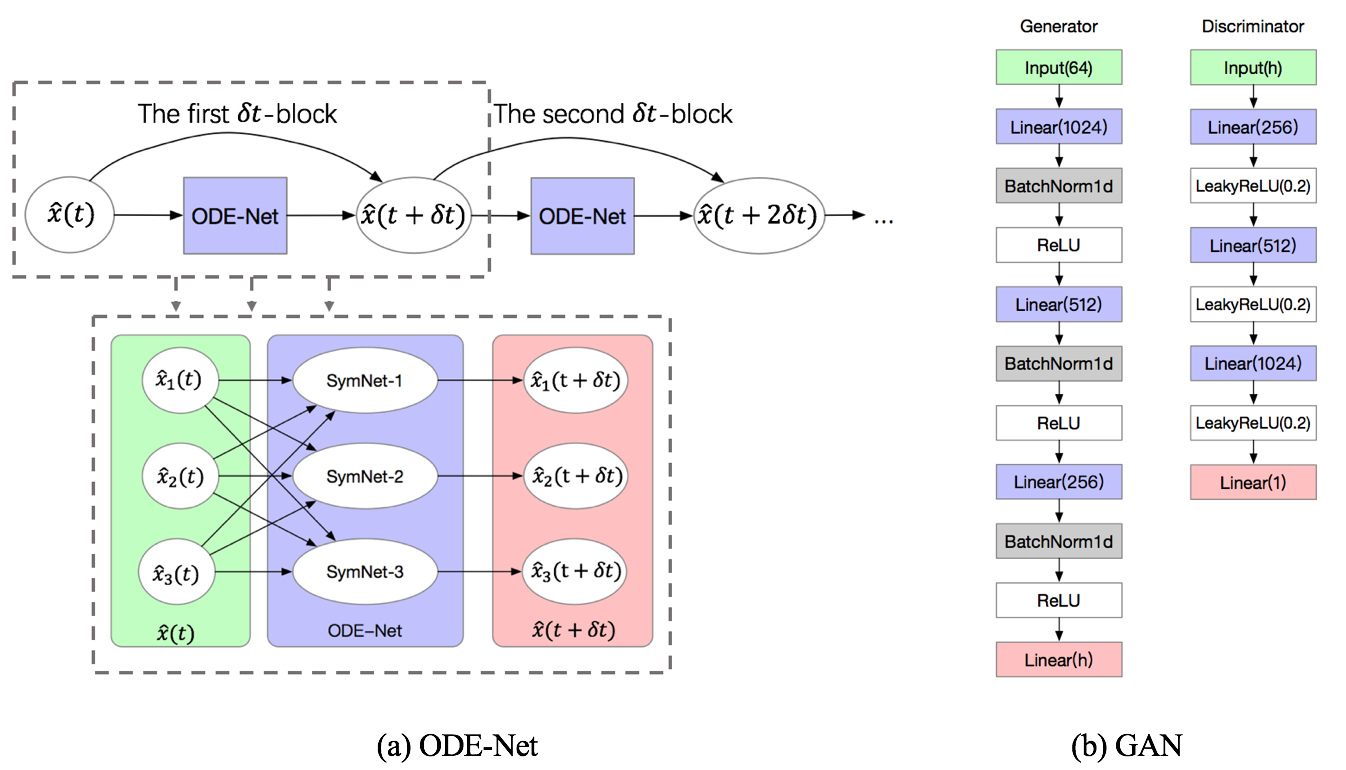}
		\caption{The architecture of the ODE-Net and the GAN}\label{fig:architecture}
\end{figure}

For each ODE-${\bm{\eta}}_i$, the associated ODE-Net is a reduced version of PDE-Net 2.0 \cite{long2019pde}, which can be constructed through 
unrolled forward Euler numerical discretization. One $\delta t$-block of the $l$-th component of the ODE-Net estimating ODE-${\bm{\eta}_i}$ can be written as:
\begin{equation}
  \hat{x}_l(t+\delta t)=\hat{x}_l(t)+\delta t\cdot\symnet^\alpha_d
  (\hat{\bm{x}}(t);\hat{{\bm{\xi}}}_i^{(l)}), 
  \hat{\bm{x}}(t)=(\hat{x}_1(t),\cdots,\hat{x}_d(t)),l=1,\cdots,d, 
\end{equation}
where $\bm{\hat{x}}(0)=\bm{x}(0)$ which is the initial value of ODE-${\bm{\eta}_i}$ and $\symnet^\alpha_d$ is a network that takes 
a $d$ dimensional vector as input and has $\alpha$ hidden layers, which can approximate polynomials of input variables and is able to output the expression of the equation through symbolic computations $s(\cdot)$. 
In practice, we take $\alpha = 2$.
Here, $\hat{{\bm{\xi}}}_i=(\hat{{\bm{\xi}}}_i^{(1)},\cdots,\hat{{\bm{\xi}}}_i^{(d)})$ are the trainable parameters of ODE-Net, and $s(\hat{{\bm{\xi}}}_i)$ is an estimation of ${\bm{\eta}}_i$. Fig. \ref{fig:architecture} (a) shows the architecture of the ODE-Net.

\subsection{The GAN Component}\label{subsec:gan}

It is natural to apply GAN to $\{s(\hat{\bm{\xi}_i})\}_i$ to estimate the distribution of $\bm{\eta}$. However, it is difficult (and unnecessary) to take derivative on the symbolic computation $s$. Therefore, in RODE-Net, we choose to use GAN to approximate the distribution of $\{\hat{\bm{\xi}_i}\}_i$ instead. By doing so, we can still provide a viable approximation of the distribution of $\bm{\eta}$ by simply taking the symbolic computations after the estimation of the distribution of $\{\hat{\bm{\xi}_i}\}_i$.

In RODE-Net, we adopt WGAN (Wasserstein GAN) \cite{goodfellow2014generative,
arjovsky2017wasserstein} to draw the distribution of ${\bm{\xi}}$ by solving 
the following constrained minimax problem: 
\begin{equation}
  \min_G \max_D \mathbb{E}_{{\bm{\xi}} \sim  \mathcal{P}_{\bm{\xi}}} D({\bm{\xi}}) - 
    \mathbb{E}_{\tilde{{\bm{\xi}}} \sim 
    \mathcal{P}_{g}} D(\tilde{{\bm{\xi}}}) 
    \quad s.t. ||D||_L \leq 1.
\end{equation}
where $D$ is a discriminator and $\mathcal{P}_g$ is the distribution implicitly defined by the generator $G$: $\tilde{{\bm{\xi}}}=G(z),z\sim p_z$. Since the true data 
${\bm{\xi}}$ are not directly observable, the training samples ${\bm{\xi}}\sim \mathcal{P}_{\bm{\xi}}$ 
are approximately generated by $\hat{{\bm{\xi}}}\sim \mathcal{P}_{\hat{{\bm{\xi}}}}$, i.e. 
the parameters of ODE-Nets which can be learned from the observed data 
(Sec\ref{subsec:ode-net}). In RODE-Net, $G$ and $D$ are neural networks with multiple fully connected 
hidden layers and non-linear activation functions. The architecture of the generator and discriminator is shown in Fig. \ref{fig:architecture} (b). The dimension of both the output of the generator and the input of the discriminator are equal to the size of each $\hat{\bm{\xi}_i}$, which is $d \cdot (\alpha^2 + 2 d \cdot \alpha +  d + 2\alpha + 1)$.

\section{Training}\label{sec:training}
The training of RODE-Net has three stages which include two warm-up stages and the alternating stage. In the first stage (warm-up-1), we train the ODE-Nets from the observed data. Secondly (warm-up-2), we train the GAN using the estimated parameters of the ODE-Nets. Thirdly (alternating stage), the GAN is used as a regularization to further improve the estimates on the parameters of the ODE-Nets by updating the parameters of GAN and ODE-Nets alternatively. 

\paragraph{Warm-up-1:}
Given each $X_i:=\{\bm{x}_{\bm{\eta}_i}^j(t_k):1\leq j\leq N_i,
0\leq k\leq S\}$ with $1\leq i\leq M$, the parameters of the associated ODE-Net can be learned by solving the following minimization problem: 
\begin{equation} \label{eq:ode-net-minimize}
  \min_{\hat{\bm{\xi}}_i}
  \frac1{N_i}\sum_{j=1}^{N_i}L_{\text{data}}(\hat{\bm{\xi}}_i, \bm{x}_{\bm{\eta}_i}^j(t)) + 
    \lambda_h L_{\text{Huber}}(\hat{\bm{\xi}}_i).
\end{equation}
Here, the data loss term 
$
  L_{\text{data}}(\hat{\bm{\xi}}_i, \bm{x}_{\bm{\eta}_i}^j(t)) = \sum_{k=1}^{\tilde S} 
    w_k \|\hat{\bm{x}}_{\bm{\eta}_i}^j(t_k)-\bm{x}_{\bm{\eta}_i}^j(t_k)||_2^2 / \delta t^2, \quad \tilde S\le S,
    $
describes the accuracy of the ODE-Net to match the data $X_i$, where $w_k=\frac{e^k}{\sum_{s=1}^{\tilde{S}} e^s}$. 
The regularization term $L_{\text{Huber}}$ is defined as follows:
\begin{equation} \label{eq:ode-net-sparsity-loss}
  L_{\text{Huber}}(\hat{\bm{\xi}}_i) = 
  \sum_{p \text{ is a component of } \hat{\bm{\xi}}_i} \ell_s(p), \quad\mbox{and}\quad 
      \ell_s(p) = \begin{cases}
    |p| - \frac{s}{2} & \text{if \,}  |p|>s \\
    \frac{1}{2s}p^2 & \text{else}
    \end{cases}.
\end{equation}
This is known as the Huber function which is a smoothed version of the absolute value. The hyperparameters $\lambda_h=s=0.001$ are used in RODE-Net.

During the training of the ODE-Nets, we adopt a strategy 
to gradually increase the number of steps $\tilde S$ from 1 to $S_{\text{max}}$. 
This strategy is to ensure the accuracy of long-term prediction and to keep the model stable during training. In practice we take $S_{\text{max}}=6$. We use Adam \cite{kingma2014adam} to minimize the above loss function with a learning rate of 0.01. 




\paragraph{Warm-up-2:} 
The learned parameters 
$\{\hat{\bm{\xi}}_i=(\hat{\bm{\xi}}_i^{(1)},\cdots,\hat{\bm{\xi}}_i^{(d)})\}$ in Warm-up-1 
serve as training data for WGAN. To improve training of WGAN, we adopt the gradient penalty methods \cite{gulrajani2017improved}, i.e. 
solving the following WGAN-GP problem which is a penalized version of the original WGAN problem: 
\begin{equation}\label{eq:wgan-gp-loss}
  \min_G \max_D \mathbb{E}_{{\bm{\xi}} \sim  \mathcal{P}_{\bm{\xi}}} D({\bm{\xi}}) - 
    \mathbb{E}_{z \sim p_z} D(G(z)) - 
    \lambda_{gp} \mathbb{E}_{{\bm{\xi}}'\sim \mathcal{P}_{{\bm{\xi}}'}}(\|\nabla_{{\bm{\xi}}'}D({\bm{\xi}}')\|_2-1)^2, 
\end{equation}
where we define $\mathcal{P}_{{\bm{\xi}}'}$ by sampling uniformly along straight lines 
between ${\bm{\xi}}\sim \mathcal{P}_{\bm{\xi}}$ and $G(z)$ with $z\sim p_z$. Moreover, $\mathcal{P}_{\bm{\xi}}$ is approximated by $\mathcal{P}_{\hat{{\bm{\xi}}}}$. We update $G$ and $D$ alternatively by optimizing the WGAN-GP loss \eqref{eq:wgan-gp-loss} where we take $\lambda_{gp}=10$.


    
    
    

\paragraph{Alternating stage:}
In the warm-up-1 stage, each of the ODE-Net is estimated only using the data $X_i$. This has two potential problems: 1) some instances of the underlying RODE may be harder to estimate than others; 2) the estimate of each ODE does not make use of the fact that it is an instance from a RODE. Therefore, to provide an estimate of the ODE instances with a uniform control of quality and to take advantage of the entire data set $X$, we use GAN as a regularization and update the parameters of GAN and the ODE-Nets alternatively. 


In addition to the data and the Huber loss term introduced in \eqref{eq:ode-net-minimize}, we introduce another loss term, GAN regularization: 
$
  L_{\text{GANreg}}(\hat{\bm{\xi}}) = D(\hat{\bm{\xi}}).
$
The regularization by GAN can be interpreted as a switching of the role of the generator and the ODE-Nets, which enforces the parameters of the ODE-Nets to be similar to what can be generated by the learned generator. The full loss used in this stage is
\begin{equation}
  L_{\text{rodenet}} (\hat{\bm{\xi}}_i; \bm{x}_{\bm{\eta}_i}(t)) = \frac1{N_i}\sum_j L_{\text{data}}(\hat{\bm{\xi}_i}, \bm{x}_{\bm{\eta}_i}^j(t)) 
  + \lambda_h L_{\text{Huber}}(\hat{\bm{\xi}}_i) + \lambda_G L_{\text{GANreg}}(\hat{\bm{\xi}}_i),
\end{equation}
where we take $\lambda_G=0.1$ in practice.
 
The alternating update of the ODE-Nets and GAN is a mixed procedure of warm-up-1 and warm-up-2, except that we have modified the loss of training ODE-Nets. The training of the RODE-Net consists of this algorithm, together with warm-up-1 and warm-up-2, and is detailed in Algorithm \ref{alg:stage2}.

\begin{algorithm}
\caption{RODE-Net Training: rode-net-train($\theta, V, \bm{\xi}, \bm{x}(t))$. \\ 
The hyperparameters are set to $n_d=1, n_{\text{GAN}}=100$. The parameters of ODE-Nets and the GAN are initialized by independent gaussian distribution. }
\label{alg:stage2}
\begin{algorithmic}[1]
\STATE \textbf{Warm-up-1}: Compute an initial estimation of the ODE-Nets' parameters $\{\hat{\bm{\xi}}_i\}$;
\STATE \textbf{Warm-up-2}: Compute an initial estimation of the GAN's parameter $\theta,V$;
\WHILE{$\{\hat{\bm{\xi}}_i\}$ has not converged}
    \FOR{$i$=1,...,$M$}
    \STATE Update the $i$-th ODE-Net using $L_{\text{rodenet}} (\hat{\bm{\xi}}_i; \bm{x}_{\bm{\eta}_i}(t))$ ;
    \ENDFOR
    
    \FOR{Iteration of discriminator $n_{d}$}
    \STATE Update the discriminator $D_V$ using WGAN-GP loss ;
    \ENDFOR
    \FOR{Iteration of GAN $n_{\text{GAN}}$}
    \STATE Update the GAN using $\{\hat{\bm{\xi}}_i\}$ as training data;
    \ENDFOR
\ENDWHILE
\end{algorithmic}
\end{algorithm}

\section{Experiments}
In this section, we show the effectiveness of the proposed RODE-Net on simulation data. We first describe how simulation data is generated. Then, we show
how well the trained RODE-Net estimates the analytical form and the parameter distribution of the unknown RODE and generates reliable predictions.

\subsection{Simulation Data}\label{simulation_data}
We select the following RODE to generate the observed data. The chosen RODE is a quadratic equation with three random parameters.
\begin{equation}
\frac{dx_1}{dt} = a(x_2-x_1),\quad
\frac{dx_2}{dt} = x_1(b-x_3)-x_2,\quad
\frac{dx_3}{dt} = x_1x_2-cx_3.
\label{quadratic_data}
\end{equation}
The observed variable $\bm{x} = [x_1, x_2, x_3]^T$ is a 3-dimensional vector. We choose two different settings to describe the randomness of the parameters $a, b, c$ in this RODE:
\begin{enumerate}
    \item \textbf{RODE\_ind:} $a \sim N(2,1), b \sim N(-1, 4), c \sim N(1, 1)$; 
    \item \textbf{RODE\_dep:} {\small
    $\begin{pmatrix} a \\ b \\ c \end{pmatrix} \sim N\left(\begin{pmatrix} 2 \\ -1 \\1 \end{pmatrix},
    \begin{bmatrix}
    1 & -2 & 1 \\
    -2 & 4 &-2 \\
    1 & -2 & 1
    \end{bmatrix}\right)$}. 
\end{enumerate}

We sample $M=500$ groups of the random parameters to form the sampled ODEs from this RODE and randomly chosen $N_i=5$ initial points from the uniform cube $U([-10,10]^3)$ for each ODE. For each initial point, we evolve the dynamics using fourth order Runge-Kutta method for $S=50$ time steps with step size $dt=0.05$. Moreover, noise is injected in the numerical solutions of these ODE instances by $\tilde{\bm{x}}(t) = \bm{x}(t) + n_r \sqrt{Var(\bm{x}(t))} w$,$w \sim N(0,I)$ with three different noise levels $n_r \in \{0, 0.01, 0.03\}$. 

In this section, prediction error at time $t$ is defined as the relative error between the predicted data $\hat{\bm{x}}(t)$, which is generated by the trained RODE-Net with input $\bm{x}(0)$, and the true data $\bm{x}(t)$:
$\mbox{Error}(t) = \frac{||\hat{\bm{x}}(t)-\bm{x}(t)||_2}{ ||\bm{x}(t) ||_2}$.

\subsection{The performance of RODE-Net}
This subsection presents the performance of the learned RODE-Net in terms of prediction and identification of the random parameters.

\begin{figure}[htbp!]
		\centering
		\includegraphics[width = 0.8\textwidth]{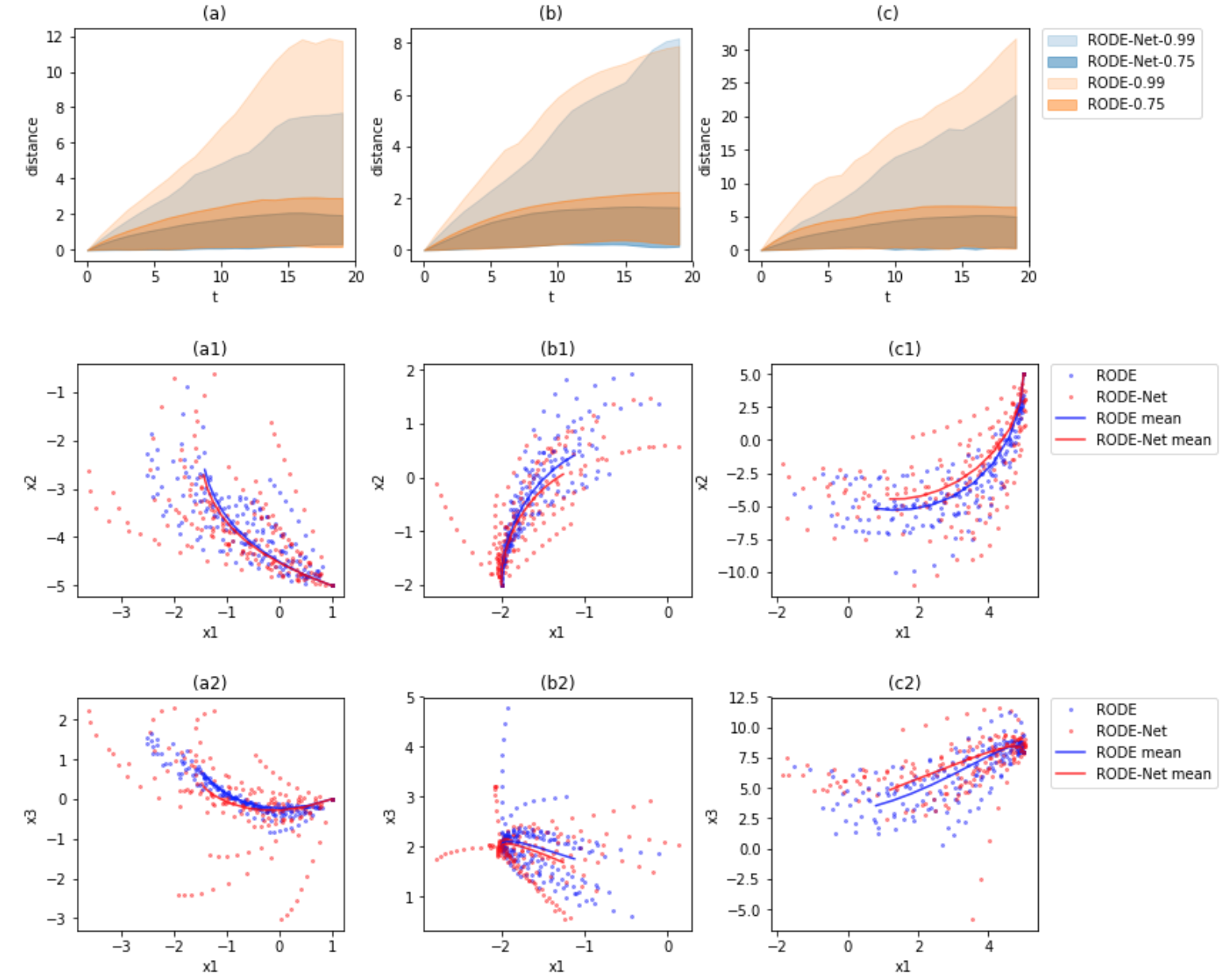}
        \caption{The predicted trajectories of the ODEs sampled from RODE-Net and the original RODE.}\label{fig:GAN_pred}
\end{figure}

To validate the performance of RODE-Net on prediction, we use the learned GAN to generate several ODEs and inspect whether they distribute similarly to those sampled from the true RODE. Fig. \ref{fig:GAN_pred} shows the predicted trajectories of 500 ODEs generated by the RODE-Net (labeled as RODE-Net) and 500 sampled instances of the RODE  \eqref{quadratic_data} (labeled as RODE) using three randomly selected initial points. 
At each time step $t$, we calculate the Euclidean distance between each trajectory and their mean and plot the 0-99\% and 0-75\% bands of the distance in Fig. \ref{fig:GAN_pred} (a-c). 
The bands of two groups generally overlap indicating that the distribution of the predictions of RODE-Net is similar to those of the RODE.
This conclusion can also be validated by the rest of the plots in Fig. \ref{fig:GAN_pred}, where we can see the data points on the trajectories of the sampled RODE-Net and RODE are distributed in similar regions.
The solid lines in Fig. \ref{fig:GAN_pred}(a1)-(c2) represent the mean of the trajectories. The closeness of the two solid lines in each of these subfigure indicates that the bias between the predictions of the RODE-Net and that of the true RODE is relatively small. 

We also compare the estimated distribution of the parameters from RODE-Net with those of the RODE by calculating the mean and standard deviation of the coefficients from 100 generated ODEs from RODE-Net. The comparisons for the case RODE\_ind are shown in Table \ref{RODE-Net-quadratic}. We can see that the RODE-Net successfully differentiates between the terms with and without randomness. For the terms exist in the RODE, the mean and standard deviation estimated by the RODE-Net is close to those of the RODE, while the estimated coefficients that do not exist in the RODE have relatively small means and standard deviations. 
\begin{table}[htbp!]
\caption{The mean and standard deviation of the estimated coefficients for RODE\_ind.}
\label{RODE-Net-quadratic}
\centering
\begin{small}
\begin{tabular}{c|c|c|c|c|c|c}     
\toprule
     \multirow{2}*{term} & \multicolumn{2}{c|}{coefficients of $\frac{dx_1}{dt}$} &  \multicolumn{2}{c|}{coefficients of $\frac{dx_2}{dt}$}&  \multicolumn{2}{c}{coefficients of $\frac{dx_3}{dt}$} \\
     \cline{2-7}
     & RODE-Net & RODE\_ind & RODE-Net & RODE\_ind & RODE-Net & RODE\_ind \\
    \hline  
    
   1 & 0.029(0.764) & 0  & -0.069(0.527) & 0 & 0.075(0.904)& 0\\
  \hline  
  $x_1$ & -1.841(0.921) & -2.000(1.000) & -0.770(1.725) &-1.000(2.000) &0.025(0.414)& 0 \\
  \hline  
  $x_2$ & 1.729(0.873) & 2.000(1.000) & -0.988(0.304) & -1.000(0) & -0.051(0.506)& 0\\
  \hline  
  $x_3$ & -0.002(0.051)& 0 & -0.007(0.140) & 0 &-0.904(0.834)& -1.000(1.000) \\
  \hline  
  $x_1x_3$ & -0.035(0.027)& 0 & -0.805(0.180) & -1.000(0) & -0.003(0.070) & 0 \\
  \hline  
  $x_1^2$ & 0.006(0.058)& 0 & 0.035(0.267)& 0 & -0.077(0.149)& 0 \\
  \hline  
  $x_1x_2$ & -0.001(0.065)& 0 & 0.020(0.198) & 0 & 0.799(0.397)& 1.000(0) \\
  \bottomrule
\end{tabular}
\end{small}
\end{table}
\begin{table}[htbp!]
\caption{The estimated covariance of the coefficients for RODE\_ind and RODE\_dep.}
\label{GAN_cov}
\centering
\begin{tabular}{c|c|c|c|c|c|c}     
\toprule
     \multirow{2}*{} & 
     \multicolumn{2}{c|}{$cov(a,b)$} &  \multicolumn{2}{c|}{$cov(a,c)$}&  \multicolumn{2}{c}{$cov(b,c)$} \\
     \cline{2-7}
     & estimated & true & estimated & true & estimated & true \\
    \hline  
   RODE\_ind & 0.075 & 0 &-0.019  & 0 & -0.061 & 0\\
  \hline  
  RODE\_dep &-1.659  & -2 & 0.900  & 1 & -1.700 & -2 \\
  \bottomrule
\end{tabular}
\end{table}

Moreover, we show that RODE-Net can also well identify RODE with a joint distribution of its parameters. Table \ref{GAN_cov} shows the covariance of the coefficients of 100 ODEs generated from RODE-Net. We use the coefficient of term $x_2$ of $\frac{dx_1}{dt}$, $x_1$ of $\frac{dx_2}{dt}$, $-x_3$ of $\frac{dx_3}{dt}$ to estimate the random coefficients $a, b, c$ in \eqref{quadratic_data}. In both cases, the covariance of the coefficients of RODE-Net are close to the true value.


\subsection{The Regularization Effect of GAN in RODE-Net}

In this subsection, we demonstrate the regularization effect of the trained GAN in improving the estimates of ODE-Nets. We compare the performance of the RODE-Net (i.e. ODE-Nets with GAN regularization trained by Algorithm \ref{alg:stage2}) with the trained ODE-Nets without GAN regularization (trained by warm-up-1 of Algorithm \ref{alg:stage2}) and find that the GAN regularization helps in reducing prediction error and correcting the learned formulas of the ODEs. For convenience, we shall refer to the ODE-Nets without GAN regularization simply as ODE-Nets.

\begin{figure}[htbp!]
		\centering
		\includegraphics[width = 0.75\textwidth]{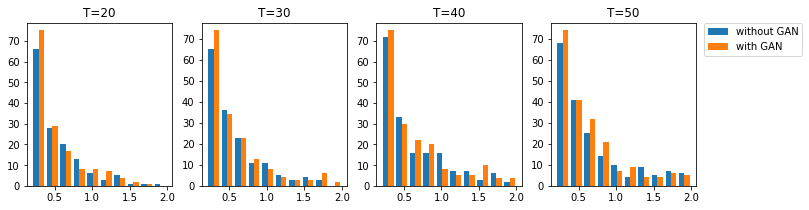}
		\caption{The histogram of the prediction error of the ODE-Nets trained with and without GAN regularization at four time instances.}\label{fig:stage2_overall_error}
\end{figure}

Let $e^i_{1}(t)$ (resp. $e^i_{2}(t)$) be the median of relative errors of the prediction of the $i$-th sample of the RODE-Net (resp. ODE-Nets) at time $t$ among 100 initial points. We generate in total $500$ instances. Denote the error vectors $\vec{e}_j(t)=(e^1_j(t),\ldots,e^{500}_j(t))$ with $j=1,2$. We compare the histograms of $\vec{e}_1(t)$ and $\vec{e}_1(t)$ at $t=20, 30, 40, 50$ in Fig. \ref{fig:stage2_overall_error}. For a better visual comparison, we cropped the vectors $\vec{e}_j(t)$ within the range $[0.2,2]$. The numbers of ODE instances for the
RODE-Net and ODE-Nets with error $<0.2$ are comparable, while the number of ODE instances for RODE-Net is smaller than that of the ODE-Nets with error $>2$. 

\begin{figure}[htbp!]
		\centering
		\includegraphics[width = 0.75\textwidth]{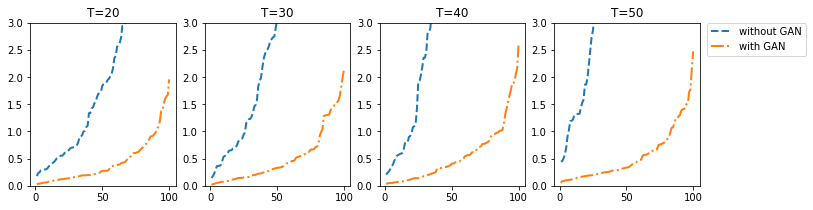}
		\caption{The order statistic of prediction error at four time instances.}\label{fig:statistic_prediction_error}
\end{figure}

As one can see that, the GAN regularization of RODE-Net notably improves overall prediction errors in comparison with ODE-Nets. For some ODE instance the improvement can be rather significant. In Fig. \ref{fig:statistic_prediction_error}, we plot the quantiled relative error vectors for a particular instance of the ODE-Nets with and without GAN regularization.


\begin{figure}[htbp!]
\centering
\begin{subfigure}{0.4\textwidth}
\includegraphics[width = \textwidth]{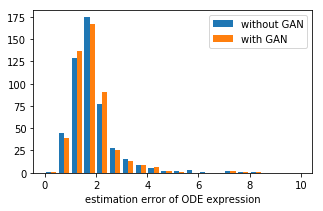}
\caption{ }\label{fig:expression_error_compare}
\end{subfigure}
$\quad$
\begin{subfigure}{0.50\textwidth}
\begin{tabular}{c|c|c|c}     
\toprule
term & without GAN & with GAN & true \\
\hline  
$x_1x_2$ & 0.9074 & 0.9164 & 1 \\
\hline  
$1$ & 0.3755 & 0.0793 & 0 \\
\hline  
$x_2$ & -0.2055 & -0.1316 & 0 \\
\hline  
$x_3$ & -0.2567 & -0.1651 & -0.143 \\
\bottomrule
\end{tabular}
\caption{ }\label{stage2_expression_boosting}
\end{subfigure}
\caption{Left: The histograms of $\vec{e}_{\bm{\eta}}$ of ODE-Nets with and without GAN regularization. the mean and standard deviation of the errors among 500 ODE instances of ODE-Nets with and without GAN regularization are 2.093(3.352) and 2.122(3.372) respectively. Right: The coefficients of $\frac{dx_3}{dt}$ in the learned ODE.}
\end{figure}

GAN regularization in RODE-Net also helps with the estimation of the expression of the RODE. In Fig. \ref{fig:expression_error_compare}, we compared the distribution of the errors $\vec{e}_{\bm{\eta}}:=\{\|\hat{\bm{\eta}_i}-\bm{\eta}_i\|_1: 1\le i\le 500\}$ associated to the ODE-Nets with and without GAN regularization. As one can see that, GAN regularization indeed helps with the estimation of the expression of the RODE. 


We note that some of the estimations by RODE-Net is notably more accurate than the ODE-Nets without GAN regularization. An example is shown in the table in Fig.\ref{stage2_expression_boosting}. We can see that the estimated coefficients of the primary terms of the ODE are closer to the true values when trained with GAN regularization.

\section{Conclusion and Future Work}
In this paper, we proposed a new framework called RODE-Net to identify RODE from data. The RODE-Net is a combination of ODE-Net for system identification and GAN for parameter distribution estimation. Numerical experiments demonstrated that the trained RODE-Net was able to identify ODEs with minor prior knowledge on the dynamics, but was also able to well estimate the distribution of the RODE. This enabled us to simulate trajectories using the trained RODE-Net that distributed similarly to the ODEs generated from the true RODE. Moreover, the well trained GAN in the RODE-Net served as an effective regularization that enabled the RODE-Net to obtain more accurate predictions and parameter estimations compared to the ODE-Nets without GAN regularization. One limitation and possible further extension of the RODE-Net is that, due to the high dimensionality of the parameters of the ODE-Net, one may include more prior knowledge to reduce the dimension. Furthermore, we would also like to apply the proposed framework to more practical examples, such as the particle fluid systems. It is also worth to study how to design a better data driven regularization based on or beyond GAN. 

\section{Broader Impact}
The proposed RODE-Net differs from earlier works major in twofold and hence it may lead to further impact. 
1) The proposed RODE-Net was used to estimate the parameters' distribution of the RODE rather than the best set of parameters in average sense. While most classical methods focus on uncertainty forward propagation and parameter calibration, the RODE-Net may encourage us to investigate alternative ways to model systems with inherent randomness. 2) The learned GAN was further used as a data driven regularization. This approach can also be generalized to inverse problems of partial differential equations. It is worth exploring on how to use GAN to improve or even replace the commonly used empirical regularization methods. However, due to the non-transparency of deep neural networks, the learned generator of GAN and the corresponding data driven regularization is not always reliable. One should be cautious when deploying GAN based regularization in practice, especially for systems with low fault tolerance. 

\begin{ack}
Zichao Long is supported by The Elite Program of Computational and Applied Mathematics for PhD Candidates of
Peking University. Bin Dong is supported in part by Beijing Natural Science Foundation (No. 180001) and Beijing Academy of Artificial Intelligence (BAAI).
\end{ack}
\small

\begin{thebibliography}{10}

\bibitem{long2019pde}
Zichao Long, Yiping Lu, and Bin Dong.
\newblock Pde-net 2.0: Learning pdes from data with a numeric-symbolic hybrid
  deep network.
\newblock {\em Journal of Computational Physics}, page 108925, 2019.

\bibitem{lin2008learning}
Zhouchen Lin, Wei Zhang, and Xiaoou Tang.
\newblock Learning partial differential equations for computer vision.
\newblock {\em Peking Univ., Chin. Univ. of Hong Kong}, 2008.

\bibitem{liu2010learning}
Risheng Liu, Zhouchen Lin, Wei Zhang, and Zhixun Su.
\newblock Learning pdes for image restoration via optimal control.
\newblock In {\em European Conference on Computer Vision}, pages 115--128.
  Springer, 2010.

\bibitem{long2018pde}
Zichao Long, Yiping Lu, Xianzhong Ma, and Bin Dong.
\newblock Pde-net: Learning pdes from data.
\newblock In {\em International Conference on Machine Learning}, pages
  3214--3222, 2018.

\bibitem{patel2018nonlinear}
Ravi~G Patel and Olivier Desjardins.
\newblock Nonlinear integro-differential operator regression with neural
  networks.
\newblock {\em arXiv preprint arXiv:1810.08552}, 2018.

\bibitem{raissi2018hidden}
Maziar Raissi and George~Em Karniadakis.
\newblock Hidden physics models: Machine learning of nonlinear partial
  differential equations.
\newblock {\em Journal of Computational Physics}, 357:125--141, 2018.

\bibitem{brunton2016discovering}
Steven~L Brunton, Joshua~L Proctor, and J~Nathan Kutz.
\newblock Discovering governing equations from data by sparse identification of
  nonlinear dynamical systems.
\newblock {\em Proceedings of the National Academy of Sciences}, page
  201517384, 2016.

\bibitem{bongard2007automated}
Josh Bongard and Hod Lipson.
\newblock Automated reverse engineering of nonlinear dynamical systems.
\newblock {\em Proceedings of the National Academy of Sciences},
  104(24):9943--9948, 2007.

\bibitem{schmidt2009distilling}
Michael Schmidt and Hod Lipson.
\newblock Distilling free-form natural laws from experimental data.
\newblock {\em science}, 324(5923):81--85, 2009.

\bibitem{xu2020dlga}
Hao Xu, Haibin Chang, and Dongxiao Zhang.
\newblock Dlga-pde: Discovery of pdes with incomplete candidate library via
  combination of deep learning and genetic algorithm.
\newblock {\em arXiv preprint arXiv:2001.07305}, 2020.

\bibitem{smith2013uncertainty}
Ralph~C Smith.
\newblock {\em Uncertainty quantification: theory, implementation, and
  applications}, volume~12.
\newblock Siam, 2013.

\bibitem{liu1986probabilistic}
Wing~Kam Liu, Ted Belytschko, and A~Mani.
\newblock Probabilistic finite elements for nonlinear structural dynamics.
\newblock {\em Computer Methods in Applied Mechanics and Engineering},
  56(1):61--81, 1986.

\bibitem{zhang2001stochastic}
Dongxiao Zhang.
\newblock {\em Stochastic methods for flow in porous media: coping with
  uncertainties}.
\newblock Elsevier, 2001.

\bibitem{ghanem2003stochastic}
Roger~G Ghanem and Pol~D Spanos.
\newblock {\em Stochastic finite elements: a spectral approach}.
\newblock Courier Corporation, 2003.

\bibitem{xiu2010numerical}
Dongbin Xiu.
\newblock {\em Numerical methods for stochastic computations: a spectral method
  approach}.
\newblock Princeton university press, 2010.

\bibitem{kennedy2001bayesian}
Marc~C Kennedy and Anthony O'Hagan.
\newblock Bayesian calibration of computer models.
\newblock {\em Journal of the Royal Statistical Society: Series B (Statistical
  Methodology)}, 63(3):425--464, 2001.

\bibitem{chen2019learning}
Xiaoli Chen, Jinqiao Duan, and George~Em Karniadakis.
\newblock Learning and meta-learning of stochastic advection-diffusion-reaction
  systems from sparse measurements.
\newblock {\em arXiv preprint arXiv:1910.09098}, 2019.

\bibitem{yang2020b}
Liu Yang, Xuhui Meng, and George~Em Karniadakis.
\newblock B-pinns: Bayesian physics-informed neural networks for forward and
  inverse pde problems with noisy data.
\newblock {\em arXiv preprint arXiv:2003.06097}, 2020.

\bibitem{yang2018physics}
Liu Yang, Dongkun Zhang, and George~Em Karniadakis.
\newblock Physics-informed generative adversarial networks for stochastic
  differential equations.
\newblock {\em SIAM Journal on Scientific Computing}, 42(1):A292--A317, 2020.

\bibitem{goodfellow2014generative}
Ian Goodfellow, Jean Pouget-Abadie, Mehdi Mirza, Bing Xu, David Warde-Farley,
  Sherjil Ozair, Aaron Courville, and Yoshua Bengio.
\newblock Generative adversarial nets.
\newblock In {\em Advances in neural information processing systems}, pages
  2672--2680, 2014.

\bibitem{arjovsky2017wasserstein}
Martin Arjovsky, Soumith Chintala, and L{\'e}on Bottou.
\newblock Wasserstein generative adversarial networks.
\newblock In {\em Proceedings of the 34th International Conference on Machine
  Learning-Volume 70}, pages 214--223, 2017.

\bibitem{gulrajani2017improved}
Ishaan Gulrajani, Faruk Ahmed, Martin Arjovsky, Vincent Dumoulin, and Aaron~C
  Courville.
\newblock Improved training of wasserstein gans.
\newblock In {\em Advances in neural information processing systems}, pages
  5767--5777, 2017.

\bibitem{bora2017compressed}
Ashish Bora, Ajil Jalal, Eric Price, and Alexandros~G Dimakis.
\newblock Compressed sensing using generative models.
\newblock In {\em Proceedings of the 34th International Conference on Machine
  Learning-Volume 70}, pages 537--546. JMLR. org, 2017.

\bibitem{shah2018solving}
Viraj Shah and Chinmay Hegde.
\newblock Solving linear inverse problems using gan priors: An algorithm with
  provable guarantees.
\newblock In {\em 2018 IEEE International Conference on Acoustics, Speech and
  Signal Processing (ICASSP)}, pages 4609--4613. IEEE, 2018.

\bibitem{ledig2017photo}
Christian Ledig, Lucas Theis, Ferenc Husz{\'a}r, Jose Caballero, Andrew
  Cunningham, Alejandro Acosta, Andrew Aitken, Alykhan Tejani, Johannes Totz,
  Zehan Wang, et~al.
\newblock Photo-realistic single image super-resolution using a generative
  adversarial network.
\newblock In {\em Proceedings of the IEEE conference on computer vision and
  pattern recognition}, pages 4681--4690, 2017.

\bibitem{han2017random}
Xiaoying Han and Peter~E Kloeden.
\newblock {\em Random ordinary differential equations and their numerical
  solution}.
\newblock Springer, 2017.

\bibitem{neckel2013random}
Tobias Neckel and Florian Rupp.
\newblock {\em Random differential equations in scientific computing}.
\newblock Walter de Gruyter, 2013.

\bibitem{kingma2014adam}
Diederik~P Kingma and Jimmy Ba.
\newblock Adam: A method for stochastic optimization.
\newblock {\em arXiv preprint arXiv:1412.6980}, 2014.

\end{thebibliography}

\end{document}